\documentclass[12pt]{amsart}

\usepackage{amsmath,color}
\usepackage{slashed}
\usepackage{amsthm}
\usepackage[backref=page]{hyperref}
\usepackage[all]{xy}
\usepackage[margin=0.8in]{geometry}
\usepackage{tikz}
\usepackage{mathrsfs}

\numberwithin{equation}{section}

\makeatletter
\def\l@section{\@tocline{1}{0pt}{0pt}{2.5pc}{}}
\def\l@subsection{\@tocline{2}{0pt}{1.5pc}{3.5pc}{}}
\def\l@subsubsection{\@tocline{3}{0pt}{3.5pc}{4.5pc}{}}
\makeatother

\newtheorem{prop}{Proposition}
\newtheorem{lemma}[prop]{Lemma}

\newtheorem{thm}[prop]{Theorem}
\newtheorem{cor}[prop]{Corollary}

\numberwithin{prop}{section}

\theoremstyle{definition}

\newcommand{\del}{\partial}
\newcommand{\delb}{\bar{\partial}}
\newcommand{\dt}{\frac{\partial}{\partial t}}
\newcommand{\brs}[1]{\left| #1 \right|}

\newcommand{\gD}{\Delta}

\newcommand{\gk}{\kappa}
\newcommand{\gl}{\lambda}

\newcommand{\gw}{\omega}
\newcommand{\ga}{\alpha}
\newcommand{\gb}{\beta}

\renewcommand{\ge}{\epsilon}
\newcommand{\N}{\nabla}

\newcommand{\til}[1]{\widetilde{#1}}

\renewcommand{\bar}[1]{\overline{#1}}

\renewcommand{\i}{\sqrt{-1}}

\newcommand{\IP}[1]{\left<#1\right>}

\newcommand{\bgb}{\bar{\beta}}

\def\rsq {]\kern -2.0pt]}
\def\lsq {[\kern -2.0pt[}

\DeclareMathOperator{\inj}{inj}
\DeclareMathOperator{\tr}{tr}

\DeclareMathOperator{\Vol}{Vol}

\begin{document}

\title[Pluriclosed flow on K\"ahler manifolds]{Pluriclosed flow on K\"ahler manifolds}

\author{Jeffrey Streets}
\address{Rowland Hall\\
         University of California, Irvine\\
         Irvine, CA 92617}
\email{\href{mailto:jstreets@uci.edu}{jstreets@uci.edu}}

\date{\today}

\begin{abstract} We prove the maximal existence conjecture for pluriclosed flow on K\"ahler backgrounds.  This leads to a priori generalizations of the Aubin-Yau and Calabi-Yau theorems, and resolves the question of connectedness of generalized K\"ahler structures on Calabi-Yau backgrounds.
\end{abstract}

\maketitle

\section{Introduction}
\subsection{Cone conjecture for pluriclosed flow on K\"ahler backgrounds}

In recent decades many geometric evolution equations have been introduced (cf. \cite{phong2018new,HCF, tosatti2022chern}) with the broad aim of extending the deep theory of K\"ahler Ricci flow/K\"ahler-Einstein metrics to more general settings in complex geometry.  Among these is the pluriclosed flow \cite{PCF}, so named because it preserves the pluriclosed condition $\del\delb \omega = 0$.  In the ensuing years it was shown to be linked to the stringy renormalization group flow \cite{PCFReg}, Hitchin's generalized geometry program \cite{StreetsPCFBI, JFS}, and lead to a geometrization conjecture for complex surfaces \cite{st-geom}.  The starting point for the analysis of K\"ahler-Ricci flow is the maximal existence time criterion of Tian-Zhang \cite{TianZhang}.  Inspired by this the author and Tian conjectured a similar behavior for the pluriclosed flow based on the positive Aeppli cone \cite[Conjecture 5.2]{PCFReg}.  Our first main result confirms this conjecture on all K\"ahler manifolds:
\begin{thm} \label{t:mainthm1} Fix $(M^{2n}, J)$ a compact K\"ahler manifold, $\gw_0$ a pluriclosed metric on $M$, and let
\begin{align*}
\tau^* = \sup \{ t > 0\ |\ [\gw_0] - t c_1(M) \in \mathcal P \}.
\end{align*}
The solution to pluriclosed flow with initial condition $\gw_0$ exists smoothly on $[0, \tau^*)$.
\end{thm}

\noindent The Tian-Zhang theorem exploits the reduction of K\"ahler-Ricci flow to a scalar parabolic complex Monge-Amp\`ere equation, building on the well-developed theory of such equations, notably the foundational works of Yau \cite{YauCC} and Cao \cite{CaoKRF}.  Roughly speaking the analysis of such equations splits into three main parts: a $C^0$ estimate for the potential via maximum principle or Yau's Moser iteration argument, a uniform ellipticity/parabolicity estimate via a Pogorelov-type maximum principle, and finally higher order regularity through Yau's $C^3$ estimate \cite{YauCC} and Schauder theory.

The proof of Theorem \ref{t:mainthm1} requires the development of a comprehensive structural extension of the  scalar Monge-Amp\`ere theory to a more general \emph{system} of equations.  Most natural classes of Hermitian, non-K\"ahler metrics do not admit local scalar potentials (a notable exception being generalized K\"ahler metrics \cite{alvarez2024symplectic,MP2}).  In particular, a pluriclosed metric is determined locally by a $(1,0)$-form, hence the pluriclosed flow is genuinely a parabolic system of equations, and cannot be reduced to a scalar Monge-Amp\`ere flow.  Due to this in \cite{StreetsPCFBI} the author introduced a reduction of the pluriclosed flow inside the positive cone to a pair $(\beta, f)$ of a $(1,0)$-form and a smooth function.  The time-dependent metric is expressed $\gw = \hat{\gw} + \delb \gb + \del \bgb + d d^c f$, where $\hat{\gw}$ is a moving background. Roughly speaking the function $f$ obeys a complex Monge-Amp\`ere type flow coupled to $\beta$, which captures the failure of the metric to be K\"ahler.

A key difficulty is the lack of a direct estimate for the potential function $f$ via either maximum principle techniques or Yau's $C^0$ estimate.  To get control we develop an extension of the classic Lelong-Jensen estimate (cf. \cite{Demailly82}) to obtain a Gaussian bound for the normalized potential (Proposition \ref{p:Gaussian}).  The extra term $\beta$ gives must be controlled when integrating against the Greens function.  Estimating this term requires a weighted bound on the derivative of the Greens function, which we achieve by generalizing recent work in K\"ahler geometry \cite{GuedjTo25, guo2024diameter}.  This yields a Gaussian-type bound for the normalized potential, weighted by the supremum of $\brs{\gb}^2_g$.  On K\"ahler backgrounds there is a direct maximum principle for this quantity (Proposition \ref{p:betabounds}), which then yields the required Gaussian estimate.

While this Gaussian bound on $f$ is quite strong, it is still an integral estimate and does not directly yield uniform parabolicity via the maximum principle.  Instead, we do obtain a priori $L^p$ estimates on the inverse eigenvalues of the metric, and improve this via an adaptation of Moser iteration, generalizing techniques developed in \cite{bamler2017heat,wang2008conditions}.  Specifically, we exploit $L^p$ norms with a weighted measure defined using the auxiliary scalar satisfying the \emph{dilaton flow} (cf. \cite{Streetsscalar}).  With this weight the time derivative of $L^p$ norms is a natural combination of the linearized operator and the time-dependent \emph{weighted scalar curvature}.  This quantity is has monotone nondecreasing infimum \cite{Streetsscalar}, generalizing the classic observation of Hamilton for Ricci flow \cite{Hamilton3folds}.  This lower bound controls the relevant nonlinear term, allowing for an adaptation of parabolic Moser iteration to yield an a priori $L^{\infty}$ bound on the eigenvalues, and hence uniform parabolicity.  With this in place we apply the higher order regularity for pluriclosed flow \cite{StreetsPCFBI}.  Note that no general theory (i.e. Evans-Krylov) applies here as again it is a parabolic \emph{system}, and the arguments of \cite{StreetsPCFBI} are rather an extension of the strategy of Yau's $C^3$ estimate \cite{YauCC}, exploiting the relationship of pluriclosed flow to generalized geometry (cf. \cite{JordanStreets, JFS}).

\subsection{Aubin-Yau and Calabi-Yau Theorems}

The first geometric applications of the K\"ahler-Ricci flow were shown by Cao \cite{CaoKRF}, who first gave a new proof of the Aubin-Yau \cite{Aubin,YauCC} Theorem.  Our next main theorem proves the convergence of the normalized pluriclosed flow to a K\"ahler-Einstein metric on compact complex manifolds with $c_1(M) < 0$, in particular reproving the Aubin-Yau Theorem.

\begin{thm} \label{t:AubinYau} Fix $(M^{2n}, J)$ a compact K\"ahler manifold with $c_1(M) < 0$.  Then there exists a unique K\"ahler-Einstein metric $\gw_{\mathrm{KE}} \in - c_1(M)$, and given $\gw_0$ any pluriclosed metric on $M$, the solution to the normalized pluriclosed flow converges to $\gw_{\mathrm{KE}}$.
\end{thm}

\noindent The proof requires refining the estimates of Theorem \ref{t:mainthm1}.  The a priori estimate for $\beta$ immediately improves to give exponential decay in the normalized flow.  Crucially we require a uniform estimate for the average of $f$ to employ the Gaussian bound, which is achieved via monotonicity for a relative entropy functional inspired by Berman \cite{berman2013relative}.  Finally, we note that the a priori scalar curvature bound worsens under the normalized flow, which makes the the Moser iteration argument as employed in Theorem \ref{t:mainthm1} degenerate at infinite time.  However, by starting the iteration at a large value of $p$ tuned to a large time of the normalized flow, we are able to again recover uniform parabolicity of the flow, after which the convergence follows from the expander entropy functional for generalized Ricci flow \cite{Streetsexpent}.

Cao \cite{CaoKRF} also gave a K\"ahler-Ricci flow proof of the Calabi-Yau Theorem \cite{YauCC}.  Our next theorem proves the convergence of the pluriclosed flow on compact K\"ahler manifolds with $c_1(M) = 0$, in particular reproving the Calabi-Yau Theorem.

\begin{thm} \label{t:CalabiYau} Fix $(M^{2n}, J)$ a compact K\"ahler manifold with $c_1(M) = 0$.  Given $\gw_0$ a pluriclosed metric on $M$, there exists a unique Ricci-flat K\"ahler metric $\omega_{\mathrm{CY}} \in [\gw_0]$, and the solution to pluriclosed flow with initial condition $\gw_0$ exists for all $t \geq 0$ and converges exponentially to $\omega_{\mathrm{CY}}$.
\end{thm}

\noindent Again the proof requires a refinement of the estimates of Theorem \ref{t:mainthm1}.  A key point is to establish a uniform lower bound on the normalized potential, which is achieved via a generalization of Yau's Moser iteration argument \cite{YauCC}, further exploiting the a priori estimate for $\brs{\gb}^2_g$.  We cannot directly obtain a uniform a priori upper bound for the average value of $f$, but by a combination of integral estimates and maximum principles we show an oscillation bound for the average (Lemma \ref{l:CYaveragebnds}).  With this we obtain a refined version of the usual Schwarz Lemma estimate for the inverse trace expressed in terms of the potential function $f$ normalized by a local spacetime average.  With this we can apply the parabolic Moser iteration argument for Theorem \ref{t:mainthm1} described above uniformly in time, finishing the proof.

\subsection{Implications for generalized K\"ahler geometry}

Theorem \ref{t:CalabiYau} has implications for generalized K\"ahler geometry \cite{GualtieriGKG}, which forms the target space for $N = (2,2)$ supersymmetric sigma models backgrounds \cite{GHR}.  In \cite{GKRF} it was shown that pluriclosed flow preserves generalized K\"ahler geometry, in which case we call it generalized K\"ahler-Ricci flow (GKRF).  The GKRF on Calabi-Yau backgrounds was studied in depth in \cite{AFSU}.  In this case it was first shown that every generalized K\"ahler structure is determined by closed spinors $\psi_1, \psi_2$.  Furthermore, by \cite[Theorem 4.13]{AFSU}, given an admissible pair of spinor classes there exists a unique canonically associate Calabi-Yau structure.

Thus a natural question is whether the GKRF with given spinor classes will always limit to this Calabi-Yau metric, which in particular would imply contractibility of GK structures with fixed spinor classes on Calabi-Yau manifolds.  As explained in \cite{AFSU} such questions are subtle due to the nonlinear structure of the space of GK structures, initially observed in \cite{AGG} via Joyce deformations of hyperK\"ahler structures, and thus linked to questions in symplectic geometry such as contractibility of Hamiltonian diffeomorphism groups.  A priori it is not even known if this space is connected, and recent work \cite{apostolov2026poisson} hints that in general this connectedness may fail, while a positive result appears in \cite{JFS}.  Global existence and convergence results for GKRF were developed in \cite{AFSU} which assume a priori deformability of the initial data to a K\"ahler structure, hence yielding contractibility of the connected component of the Calabi-Yau structure \cite[Corollary 7.7]{AFSU}.  Theorem \ref{t:CalabiYau} immediately yields the strongest possible statement in this direction for Calabi-Yau backgrounds, namely contractibility of the space of all generalized K\"ahler structures with fixed background Poisson tensor and spinor classes.

\begin{cor} \label{c:GKcor} Fix $(M^{2n}, J)$ a compact K\"ahler Calabi-Yau manifold, $\sigma$ a holomorphic Poisson tensor, and admissible spinor classes $\ga, \gb$.  Then the generalized K\"ahler class associated to this data is contractible, i.e. $\mathcal {GK}^{J, \sigma}_{\ga,\gb} \cong \star$.

\end{cor}

\vskip 0.1in
\noindent \textbf{Acknowledgments:} The author gratefully acknowledges support from the Simons Foundation through its Travel Support for Mathematicians program, and from the NSF via DMS-2603282.  This work owes a long debt to collaborations and discussions with Vestislav Apostolov, Mario Garcia-Fernandez, Marco Gualtieri, Martin Ro\v{c}ek, Mark Stern, Gang Tian, and Micah Warren, whom we gratefully thank.

\vskip 0.1in
\noindent \textbf{AI usage statement:} The author developed the argument for Proposition \ref{p:Gaussian} through collaborative discussions with ChatGPT 6-Astra, which was also used to proofread this document.  This article contains no AI generated text and the author takes responsibility for the contents of this paper.

\section{Maximal existence on K\"ahler backgrounds} \label{s:Kahlercone}
\subsection{Moving K\"ahler background and \texorpdfstring{$(\beta, f)$}{(beta, f)}-reduction} \label{ss:reduction}

We aim to construct a reduced pluriclosed flow against a moving K\"ahler background.  First we show that on a K\"ahler manifold the K\"ahler cone is canonically identified with the positive Aeppli cone, using \cite{li2009comparing} (cf. also \cite{DemaillyPaun}).

\begin{lemma} \label{l:cones} Given $(M, J)$ a compact K\"ahler manifold, the natural map $i : H^{1,1}_{BC}(M, \mathbb R) \to H_A^{1,1}(M, \mathbb R)$ identifies the K\"ahler cone with the positive Aeppli cone.
\end{lemma}

\begin{proof} We first observe that the map is an isomorphism by the $\del \delb$-lemma.  Again using the $\del \delb$-lemma, if $\omega > 0$ is a pluriclosed metric, there exists a closed real $\gk \in \Lambda^{1,1}_{\mathbb R}$ and $\beta \in \Lambda^{1,0}$ such that
\begin{align*}
    \omega = \gk + \delb \beta + \del \bar{\beta}.
\end{align*}
It follows that
\begin{align*}
F = \omega + \del \beta + \delb \bar{\beta} = \gk + d (\beta + \bar{\beta})
\end{align*}
is closed and tames $J$.  Since $[F] = [\gk]$, \cite[Proposition 3.1]{li2009comparing} implies that $[\gk]$ contains a K\"ahler form.  It is immediate that the K\"ahler cone is contained in the positive Aeppli cone, finishing the proof.
\end{proof}

Using Lemma \ref{l:cones} we can now set up pluriclosed flow solutions against a moving K\"ahler background.  Fix $\tau < \tau^*$, and fix a K\"ahler metric $\omega_K \in [\omega_0]$.  For any background Hermitian metric $h$, it follows that there exists $\chi_{\tau} \in C^{\infty}(M)$ such that
\begin{align*}
    \omega_K - \tau \rho(h) + d d^c \chi_{\tau} > 0.
\end{align*}
Then define
\begin{align*}
    \hat{\gw}_t = \gw_K - t \rho(h) + d d^c \chi_t, \qquad \chi_t = \frac{t}{\tau} \chi_{\tau}.
\end{align*}
It follows that $\hat{\gw}_t$ is a one-parameter family of K\"ahler metrics such that $\hat{\gw}_t \in [\omega_0] - t c_1(M)$, and so serves as a background against which to solve pluriclosed flow.

For any smooth existence interval of the pluriclosed flow $\omega_t$, we solve for a family of pairs $(\beta, f)$ of $(1,0)$-forms and smooth functions,
\begin{gather} \label{f:reducedflow}
\begin{split}
    \square \beta =&\ - T_{g_t} \circ \delb \beta,\\
    \dt f =&\ \log \frac{\gw_t^n}{\gw_h^n} - \tr_{\gw_t} (\del \bar{\gb} + \delb \beta) - \dot{\chi},
\end{split}
\end{gather}
where $T_{g_t}$ is the time-dependent Chern torsion and 
\begin{gather*}
    \square = \dt - \Delta^C_{g_t}
\end{gather*}
is the heat operator associated to the time-dependent Chern Laplacian $\Delta^C_{g_t}$, suitably interpreted on smooth functions and all form bundles.  These equations are the refinement of the $1$-form reduction of pluriclosed flow introduced in \cite{StreetsPCFBI}, with two key improvements coming from the K\"ahler background.  First, the moving background has no torsion, and second, through the $\del\delb$-lemma we moved the forcing term $\dot{\chi}$ to the scalar potential equation, whereas generically it appears in the $\beta$ evolution (i.e. the term $\mu$ in the notation of \cite{StreetsPCFBI}).  In particular, following the computations of \cite[Proposition 3.9]{StreetsPCFBI} we have 
\begin{align} \label{f:flowreconstruction}
    \gw_t = \hat{\gw}_t + \delb \beta_t + \del \bar{\beta}_t + d d^c f_t,
\end{align}
for all smooth existence times $t$.

We end this subsection with a lemma which estimates the volume in terms of the $\beta$-potential in a spirit similar to Lemma \ref{l:cones}.

\begin{lemma} \label{l:volumebound} Fix $(M^{2n}, \hat{g}, J)$ a compact K\"ahler manifold, and $g$ a pluriclosed metric on $M$ satisfying
\begin{align*}
    \gw =&\ \hat{\gw} + \delb \beta + \del \bar{\beta} + d d^c f
\end{align*}
Then $c \Vol(\hat{g}) \leq \Vol(g) \leq \Vol(\hat{g})$, where $c^{-1} = C_n(1 + \sup_M \brs{\del \beta}_g^n)$.
\begin{proof}
    As in Lemma \ref{l:cones} we define the closed form $F = \gw + \del \beta + \delb \bgb = \hat{\gw} + d(\beta + \bar{\beta}) + d d^c f$.  An exterior algebra computation (cf. \cite[Prop. 5.10]{PCF}) taking note of the $(2,0)-(0,2)$ components of $F$ shows
    \begin{align*}
        \frac{F^n}{n!} = \Phi dV_g, \qquad \Phi = \sum_{j=0}^{\lfloor{n/2}\rfloor} \frac{\brs{(\del \beta)^{\wedge j}}^2_g}{(j!)^2}.
    \end{align*}
Noting that $\int_M F^n = \int_M \hat{\gw}^n$ by Stokes Theorem and using interpolation on $\Phi$ gives the lemma.
\end{proof}
\end{lemma}

\subsection{Maximum principle and Harnack estimates}

First we record an a priori lower bound for $f$ which follows directly from the maximum principle.

\begin{prop} \label{p:flower} Given the setup above, one has
\begin{align*} 
    \square f =&\ \tr_{g_t} \hat{g} + \log \frac{\gw_t^n}{\gw_h^n} - n - \dot{\chi}.
\end{align*}
Moreover, for any smooth existence time $t$,
\begin{align*}
    \inf_{M \times \{t\}} f \geq \inf_{M \times \{0\}} f - C t.
\end{align*}
\begin{proof} The claimed heat equation follows by combining (\ref{f:reducedflow}) and (\ref{f:flowreconstruction}).  If we let $\gl_i$ denote the eigenvalues of $g^{-1} \hat{g}$ we then obtain
\begin{align*}
    \square f = \sum_i \left(\gl_i - 1 - \log \gl_i \right) + \log \frac{\hat{\gw}^n}{\gw_h^n} - \dot{\chi} \geq - C,
\end{align*} 
and the estimate follows by the maximum principle.
\end{proof}
\end{prop}

Next we record an a priori estimate for $\beta$ and $\del \beta$ (cf. \cite[Lemma 6.2]{AFSU}).

\begin{prop} \label{p:betabounds} Given the setup above, one has
\begin{align*}
    \square \brs{\beta}^2_g =&\ - \brs{\N \beta}^2_g - \brs{\bar{\N} \beta + \beta \cdot \bar{T}}^2_g\\
    \square \brs{\del \beta}^2_g =&\ - \brs{\N \del \beta}^2_g - \tfrac{1}{2} \brs{T}^2_g - \IP{Q, \tr_g \del \beta \otimes \delb \bar{\beta}}_g.
\end{align*}
In particular, for any smooth existence time $t$,
\begin{align*}
    \sup_{M \times \{t\}} \brs{\beta}^2_g \leq&\ \sup_{M \times \{0\}} \brs{\beta}^2_g\\
    \sup_{M \times \{t\}} \brs{\del \beta}^2_g \leq&\ \sup_{M \times \{0\}} \brs{\del \beta}^2_g.
\end{align*}
\end{prop}

We also record the scalar curvature monotonicity for generalized Ricci flow observed in \cite{Streetsscalar}.  For this we need an auxiliary function, the dilaton flow, defined here as the unique solution to
\begin{align} \label{f:dilatonflow}
    \square \phi = \tfrac{1}{2} \brs{T}^2_g, \qquad \phi_0 = 0.
\end{align}
In the notation of \cite{Streetsscalar}, and accounting for the time scaling by $1/2$, the right hand side would be expressed as $\tfrac{1}{12} \brs{H}^2_g$, where $H$ is the Bismut torsion, but this expression is more convenient for us here.  For the dilaton weight $\phi$ we obtain the weighted scalar curvature
\begin{align} \label{f:scalarcurvature}
    R^{H,\phi} = R - \tfrac{1}{12} \brs{H}^2 + 2 \Delta_g \phi - \brs{\N \phi}^2.
\end{align}
We further observe the key formula arising from \cite[Lemma 3.2(1)]{Streetsscalar}, taking into account gauge-fixing by the Lee vector field,
\begin{align} \label{f:phimeasureevolution}
    \frac{d}{dt} \int_M u e^{-\phi} dV_g = \int_M \left( \square u - \tfrac{1}{2} R^{H,\phi} u \right) e^{-\phi} dV_g.
\end{align}
We next record a priori estimates for these objects.

\begin{prop} \label{p:dilatonscalarbound} Given the setup above, one has for all smooth existence times $t$,
\begin{align*}
    0 \leq \phi \leq C, \qquad R^{H,\phi} \geq \inf_{M \times \{0\}} R^{H,\phi}.
\end{align*}
\begin{proof} The scalar curvature lower bound is \cite[Corollary 2.5]{Streetsscalar}.  The lower bound for $\phi$ follows immediately from the maximum principle.  For the upper bound we exploit Proposition \ref{p:betabounds} to observe
\begin{align*}
    \square \left( \phi + \brs{\del \beta}^2_g \right) \leq 0,
\end{align*}
and apply the maximum principle.
\end{proof}
\end{prop}

\subsection{A coupled Gaussian bound}

\begin{prop} \label{p:Gaussian} Fix $(M^{2n}, J, h)$ a compact K\"ahler manifold.  Suppose $\hat{\omega}$ is a background metric and fix $\beta, f$ such that
\begin{align} \label{f:metricdecomp}
    \omega = \hat{\omega} + \delb \beta + \del \bar{\beta} + d d^c f > 0, \qquad d d^c \omega = 0.
\end{align}
There are constants $c, C > 0$ depending on $h$ and $\hat{\omega}$ such that
\begin{enumerate}
    \item $\brs{\brs{f - \bar{f}_h}}_{L^1(h)} \leq C \left(1 + \sqrt{ \sup_M \brs{\beta}_g^2} \right)$,
    \item $\int_M \exp \left( \frac{c \left( f - \bar{f}_h\right)_+^2}{1 + \sup_M \brs{\gb}^2_g} \right) dV_h \leq C$,
    \item For any $\gl \geq 0$, $\int_M \exp \left( \gl (f - \bar{f}_h) \right) dV_h \leq C \exp \left( C \gl \left(1 + \gl \sup_M \brs{\gb}^2_g \right) \right)$.
\end{enumerate}
\end{prop}

This estimate is inspired the K\"ahler potential estimates obtained via Greens functions from \cite{YauCC}.  A key preliminary point is a weighted estimate for the derivative of the Greens function of $h$, which generalizes recent work in the K\"ahler setting \cite[Lemma 5.6]{guo2024diameter}, \cite[Lemma 2.5]{GuedjTo25}:

\begin{lemma} \label{l:weightedGreens} There is a constant $C$ depending on $h$ such that every pluriclosed metric $g$ satisfies
\begin{align*}
    \sup_{x \in M} \int_M \frac{ \brs{\N_h^{1,0} G^h(x,y)}_g^2}{G^h(x,y)^{n/(n-1)}}\, dV_h(y) \leq C \int_M \tr_h g\, dV_h.
\end{align*}
\begin{proof} In this proof we drop the decoration $h$ from the Greens function for clarity.  Fix some $0 < \rho_0 < \min \{1, \tfrac{1}{4} \inj(h) \}$, and further fix $0 < \ge \leq \rho_0$.  Define $\rho_x(y) = d_h(x, y)$.  Furthermore fix a cutoff function $\chi$, equal to $1$ on $[0, \rho_0^2]$ and supported in $[0, 4 \rho_0^2]$.  Now define
\begin{align*}
    \psi_{x, \ge} = \chi(\rho_x^2) \log (\rho_x^2 + \ge^2), \qquad \Theta_{x,\ge} = A \omega_h + d d^c \psi_{x,\ge},
\end{align*}
where $A$ is a constant to be determined next.  To estimate $\Theta_{x,\ge}$ we first observe via Hessian comparison for $h$ that
\begin{align*}
    d d^c \log \left( \rho_x^2 + \ge^2 \right) = \frac{2}{\rho_x^2 + \ge^2} \left( \omega_h - 2 \i \del \rho_x \wedge \delb \rho_x \right) + \frac{4 \ge^2}{(\rho_x^2 + \ge^2)^2} \i \del \rho_x \wedge \delb \rho_x + O(1).
\end{align*}
Since the cutoff function $\chi$ has bounded derivatives as well, it follows that we may choose $A$ sufficiently large so that $\Theta_{x,\ge} \geq \omega_h$ globally.  Moreover, as $\ge \to 0$ the forms $\Theta_{x,\ge}$ converge smoothly away from $x$ to a form $\Theta_x$ which moreover satisfies
\begin{align} \label{f:wG10}
    \Theta_x \geq \omega_h + 2 \rho_x^{-2} \left( \omega_h - 2 \i \del \rho_x \wedge \delb \rho_x \right), \qquad 0 < \rho_x < \rho_0.
\end{align}
The form $\Theta_x$ will provide the necessary estimate on the Greens function.  Note that the right hand side of the above estimate has eigenvalue $1$ in the radial direction, and all other eigenvalues are $1 + 2 \rho_x^{-2}$.  

First recall the Greens function estimates for $0 < \rho_x < \rho_0$ (cf. \cite[Chapter 4]{aubin1998some})
\begin{align*}
    G(x,y) \geq c \rho_x^{2 - 2n}, \qquad \delb_y G(x,y) = (2 - 2n) c_n \rho_x^{1 - 2n} \delb \rho_x + O(\rho_x^{2 - 2n}).
\end{align*}
Using this and the inequality (\ref{f:wG10}) we obtain, for $0 < \rho_x < \rho_0$,
\begin{align} \label{f:wG20}
    \frac{\brs{\N_h^{1,0} G}^2_g}{G^{n/(n-1)}} dV_h \leq C \left( \rho_x^{2 - 2n} \brs{\N^h \rho_x}_g^2 + \rho_x^{4 - 2n} \tr_h g \right) dV_h.
    \leq&\ C \omega \wedge \Theta_x^{n-1}.
\end{align}
This same inequality holds for a different $C$ on all of $M$ since for $\rho_x \geq \rho_0$ the left hand side is bounded and $\Theta_x \geq \omega_h$.  We next observe by Stokes Theorem that
\begin{align*}
    \int_M \omega \wedge \Theta_{x,\ge}^{n-1} = A^{n-1} \int_M \omega \wedge \omega_h^{n-1} = A^{n-1} \int_M \tr_h g dV_h.
\end{align*}
Thus integrating (\ref{f:wG20}), sending $\ge \to 0$, and applying Fatou's Lemma finishes the proof.
\end{proof}
\end{lemma}

\begin{proof}[Proof of Proposition \ref{p:Gaussian}] For convenience we set $B = \sup_M \brs{\gb}_g^2$.  In particular this yields $\i \gb \wedge \bar{\gb} \leq B \gw$.  Since $h$ is K\"ahler we have by Stokes Theorem
\begin{align} \label{f:Gauss5}
    \int_M \tr_h \gw dV_h = \int_M \gw \wedge \gw_h^{n-1} = \int_M \hat{\gw} \wedge \gw_h^{n-1} = \int_M \tr_h \hat{\gw} dV_h \leq C.
\end{align}
It follows easily then that
\begin{align} \label{f:Gauss3}
    \brs{\brs{\gb}}_{L^2(h)}^2 \leq C B, \qquad \brs{\brs{\gb}}_{L^1(h)} \leq C \sqrt{B}.
\end{align}
We now trace (\ref{f:metricdecomp}) with respect to $\omega_h$ and rearrange to yield
\begin{align*}
    \gD_h f = \tr_{\omega_h} \omega - \tr_{\omega_h} \hat{\omega} + 2 \tr_{\omega_h} \delb \gb.
\end{align*}
Choose now $G$ a positive Greens function for $\Delta_h$.  Using Greens formula and integration by parts then yields
\begin{align} \label{f:Gauss10}
    f(x) - \bar{f}_h = \int_M G(x,y) \left( \tr_{\gw_h} \hat{\gw}(y) - \tr_{\gw_h} \gw(y) \right) dV_h(y) + 2 \mathrm{Im} Z(x),
\end{align}
where
\begin{align*}
    Z(x) := \int_M \gb \left( \N_h^{1,0} G(x,y) \right) dV_h(y).
\end{align*}
To prove item (1), we integrate (\ref{f:Gauss10}) in $x$ and apply Fubini's theorem, standard Greens function estimates, and (\ref{f:Gauss5}).

We now turn to item (2).  We first observe using $G \geq 0$ that the term from (\ref{f:Gauss10}) ivolving $\tr_{\omega_h} \omega$ can be dropped, yielding after using integrability of $G$ that
\begin{align} \label{f:Gauss15}
    \left(f - \bar{f}_h \right)_+ \leq C + 2 \brs{Z}, \qquad \frac{\left(f - \bar{f}_h \right)_+^2}{1 + B} \leq C + \frac{8 \brs{Z}^2}{B}.
\end{align}
Note that item (2) now follows immediately if $B = 0$, so we assume $B > 0$.  First observe the pointwise inequality
\begin{align} \label{f:Gauss20}
    \brs{\gb \left( \N_h^{1,0} G \right)}^2 \leq \brs{\gb}^2_g \brs{\N_h^{1,0} G}_g^2 \leq B \brs{\N_h^{1,0} G}_g^2.
\end{align}
We now decompose the integral defining $Z$ into regions $U_{<\ge} := \{ \rho_x < \ge \}, U_{\geq \ge} := \{ \rho_x \geq \ge\}$, and call the corresponding integrals $Z_{<\ge}, Z_{\geq \ge}$.  Using (\ref{f:Gauss20}), H\"older's inequality, Lemma \ref{l:weightedGreens}, and (\ref{f:Gauss5}) we obtain
\begin{gather} \label{f:Gauss22}
    \begin{split}
        \brs{Z_{\geq \ge}(x)} \leq&\ \sqrt{B} \left( \int_M \frac{\brs{\N_h^{1,0} G}^2_g}{G^{n/(n-1)}} \right)^{\tfrac{1}{2}} \left( \int_{\rho_x \geq \ge} G^{n/(n-1)} dV_h \right)^{\tfrac{1}{2}}\\
        \leq&\ C \sqrt{B} \left( 1 + \log \frac{\rho_0}{\ge} \right)^{\tfrac{1}{2}},
    \end{split}
\end{gather}
where the last line follows by using the standard Greens function estimates to show the integrand is bounded by $C(1 + \rho_x^{-2n})$.  On the other hand we can integrate in $x$ to obtain
\begin{gather} \label{f:Gauss25}
    \begin{split}
        \brs{\brs{Z_{< \ge}}}_{L^1(h)} \leq C \int_M \brs{\gb(y)}_h \int_{d_h(x,y) < \ge} d_h(x,y)^{1-2n} dV_h(x) dV_h(y) \leq C \ge \sqrt{B},
    \end{split}
\end{gather}
using standard Greens function estimates and (\ref{f:Gauss3}).

We now fix constants to finish the proof.  Choose $a > 0$ small and $L > 0$ large, fix $\gl \geq L \sqrt{B}$, and let $\ge = \rho_0 \exp(- a \gl^2/B)$.  With these choices (\ref{f:Gauss22}) yields
\begin{align*}
    \brs{\brs{Z_{\geq \ge}}}_{L^{\infty}} \leq C \sqrt{B} \left(1 + a \gl^2 / B \right)^{\tfrac{1}{2}} \leq C \sqrt{B} +  C \sqrt{a} \gl. \leq C \left( \frac{1}{L} + \sqrt{a} \right) \gl \leq \frac{\gl}{2},
\end{align*}
where the final inequality follows by choosing $L$ and $a$ relative to background constants.  Using this fact we return to (\ref{f:Gauss25}) and use Chebyshev's inequality to now yield
\begin{align*}
    \Vol_h \{ \brs{Z} > \gl \} \leq \Vol \{ \brs{Z_{< \ge}} > \tfrac{\gl}{2} \} \leq \frac{2}{\gl} \brs{\brs{Z_{< \ge}}}_{L^1(h)} \leq C \frac{\sqrt{B}}{\gl} e^{- a \gl^2/B} \leq C e^{- a \gl^2/B}.
\end{align*}
Observe also that the same inequality will hold for $0 \leq \gl < L \sqrt{B}$ by changing $C$.  We now choose any $0 < c < a$ and apply the layer-cake formula to obtain
\begin{gather} \label{f:Gauss30}
\begin{split}
    \int_M e^{c \brs{Z}^2/B} dV_h =&\ \Vol_h(M) + \frac{2C}{B} \int_0^{\infty} \gl e^{c \gl^2/B} \Vol_h \{\brs{Z} > \gl \} d \gl\\
    \leq&\ \Vol_h(M) + \frac{2 c C}{B} \int_0^{\infty} \gl e^{- (a - c) \gl^2/B} d \gl \leq C.
\end{split}
\end{gather}
Using the estimate (\ref{f:Gauss15}) finishes the proof of item (2).

For item (3), we first observe that for $\gl \geq 0$, returning to (\ref{f:Gauss15}) we have
\begin{align*}
    \gl (f - \bar{f}_h) \leq \gl (C + 2 \brs{Z}) + C \gl + \frac{c \brs{Z}^2}{B} + \frac{\gl^2 B}{c}.
\end{align*}
Integrating and applying (\ref{f:Gauss30}) yields (3).
\end{proof}

\subsection{\texorpdfstring{$L^1$}{L1} potential estimate}

Next we estimate the background $L^1$ norm of $f$.  To do this we exploit a new moving measure which has a favorable evolution equation which exploits the a priori scalar curvature lower bound.  In fact these estimates do not require a K\"ahler background.

\begin{lemma} \label{l:scalaridentity} Given $(M^{2n}, J)$ a complex manifold, $g$ a pluriclosed metric and $\phi \in C^{\infty}(M)$, one has
\begin{align*}
    R^{H,\phi} = 2 \left( s_C + 2 \Delta^C \phi - \brs{\theta - \del \phi}^2 \right).
\end{align*}
\begin{proof} A direct computation in coordinates using identities of \cite{Ivanovstring} yields
\begin{align*}
    R = 2 s_C + \tfrac{1}{12} \brs{H}^2 - 2 \brs{\theta}^2.
\end{align*}
Furthermore there is the standard identity for Hermitian manifolds 
\begin{align*}
    \Delta_g \phi = 2 \gD^C_g \phi + 2 \mathrm{Re} \IP{\theta, \del \phi}.
\end{align*}
Combining these with the definition of $R^{H,\phi}$ yields the lemma.
\end{proof}
\end{lemma}

\begin{lemma} \label{l:weightbounds} Fix $(M^{2n}, J, h)$ a compact Hermitian manifold, and suppose $g_t$ is a pluriclosed flow on $M$, and $\phi_t$ is an associated dilaton flow.  Define a measure $\mu$ and weight $w$ via
\begin{align*}
    d \mu = e^{-2 \phi} dV_g = w dV_h.
\end{align*}
Then
\begin{align*}
\frac{\del}{\del t} \log w =&\ - \tfrac{1}{2} R^{H,\phi} - \brs{\theta - \del \phi}^2_g,\\
\frac{d}{dt} \mu(M) =&\ - \int_M \left( \tfrac{1}{2} R^{H,\phi} + \tfrac{1}{2} \brs{T}^2_g + \brs{\del \phi}^2_g \right) d \mu.
\end{align*}
Furthermore, there exist constants $C = C(g_0), C_{\tau} = C(g_0, \tau)$ such that
\begin{align*}
    w \leq C e^{C\tau}, \qquad \int_0^{\tau} \int_M \left( \tfrac{1}{2} \brs{T}^2 + \brs{\del \phi}^2_g \right) d \mu dt \leq C_{\tau}.
\end{align*}
\begin{proof} Using standard formulas for pluriclosed flow and Lemma \ref{l:scalaridentity} we obtain
\begin{align*}
    \dt \log w =&\ -2 \dt \phi + \dt \log \frac{\gw^n}{\gw_h^n}\\
    =&\ - 2 \left( \gD^C \phi + \tfrac{1}{2} \brs{T}^2_g \right) - s_C + \brs{T}^2_g\\
    =&\ - s_C - 2 \gD^C \phi\\
    =&\ - \tfrac{1}{2} R^{H,\phi} - \brs{\theta - \del \phi}^2_g,
\end{align*}
as claimed.  The evolution equation for $\mu(M)$ follows directly from this and the integration by parts identity for pluriclosed metrics
\begin{align*}
    \int_M \brs{\theta - \del \phi}^2 d\mu = \int_M \left( \tfrac{1}{2} \brs{T}^2_g + \brs{\del \phi}^2_g \right) d \mu.
\end{align*}
The claimed upper bound for $w$ follows directly from the evolution equation for $\log w$ and the a priori lower bound for $R^{H,\phi}$ from Proposition \ref{p:dilatonscalarbound}.  Since the mass $\mu(M)$ is always nonnegative, we conclude a lower bound on the time integral of $\frac{d}{dt} \mu$, whence the claimed spacetime integral bound follows from out evolution equation and the a priori lower bound for $R^{H,\phi}$.
\end{proof}
\end{lemma}

\begin{prop} \label{p:fmean} Given the setup above, one has
\begin{align*}
\bar{f}_h \leq C_{\tau}.
\end{align*}
\begin{proof} Using the a priori lower bound on $f$ from Proposition \ref{p:flower}, we define
\begin{align*}
    \til{f} = f + Ct + 1 \geq 1, \qquad A(t) = \int_M \til{f} d \mu.
\end{align*}
Now note that the a priori upper bounds for $w$ and $\phi$ from Lemma \ref{l:weightbounds} and Proposition \ref{p:dilatonscalarbound} yield
\begin{align*}
    \log \frac{\gw^n}{\gw_h^n} \leq C_{\tau}.
\end{align*}
Using this, the evolution equation (\ref{f:reducedflow}) for $f$, the evolution for $w$ from Lemma \ref{l:weightbounds}, the lower bound on $R^{H,\phi}$, the a priori estimate for $\brs{\gb}_g^2$ from Proposition \ref{p:betabounds}, and integration by parts then yields
\begin{align*}
    \frac{d}{dt} A \leq&\ C_{\tau} (1 + A) + 2 \brs{ \int_M \tr_{\gw} \left( \del \bgb + \delb \gb \right) d \mu} \\
    \leq&\ C_{\tau} (1 + A) + C_{\tau} \int_M \left( \brs{T}^2 + \brs{\del \phi}^2 \right) d \mu.
\end{align*}
Using this together with the spacetime integral estimate of Lemma \ref{l:weightbounds} yields $A(t) \leq C_{\tau}$.  We thus obtain using Propositions \ref{p:betabounds} and \ref{p:Gaussian} and the upper bound for $w$ from Lemma \ref{l:weightbounds},
\begin{align*}
    \bar{f}_h \mu(M) =&\ \int_M \left( f + \bar{f}_h - f \right) d \mu \leq C_{\tau} + C_{\tau} \brs{f - \bar{f}_h}_{L^1(h)} \leq C_{\tau}.
\end{align*}
To finish the proof it remains to obtain an a priori lower bound for $\mu(M)$, which follows from the upper bound on the dilaton from Proposition \ref{p:dilatonscalarbound} and the lower bound on total volume from Lemma \ref{l:volumebound}.
\end{proof}
\end{prop}

\begin{cor} \label{c:Lptrace} Given the setup above, for any $p \geq 1$ we have
\begin{align*}
    \sup_{t \leq \tau} \int_M (\tr_g h)^p dV_h \leq C_{\tau,p}.
\end{align*}
\begin{proof} The Schwarz Lemma argument of (\cite[Lemma 6.2]{StreetsPCFBI}, \cite[Proposition 4.12]{StreetsND}) yields
\begin{align*}
    \square \log \tr_g h \leq C \tr_g h.
\end{align*}
Using Proposition \ref{p:flower} we obtain for large $A$,
\begin{align*}
    \square \left( \log \tr_g h - A f \right) \leq - \tr_g h + C.
\end{align*}
Using again Proposition \ref{p:flower} and the maximum principle we obtain $\tr_g h \leq C_{\tau} e^{Af}$.  The claim then follows from Proposition \ref{p:Gaussian} item (3) and Proposition \ref{p:fmean}.
\end{proof}
\end{cor}

\subsection{Moser iteration and main proof}

We finish the proof of uniform ellipticity via a Moser iteration argument using a weighted measure $dm = e^{-\phi} dV_g$, and exploiting the a priori $L^p$ estimates on eigenvalues of Corollary \ref{c:Lptrace}.  The first step is a weighted Sobolev inequality.

\begin{lemma} \label{l:Sobolev} Fix $(M^{2n}, J, h)$ a compact K\"ahler manifold.  Fix $g$ a Hermitian metric on $M$ and $\phi$ a smooth function satisfying
\begin{enumerate}
    \item $\int_M \left[ (\tr_g h)^{4n^2} + \left( \frac{\gw^n}{\gw_h^n} \right)^{4n} \right] dV_h \leq A$,
    \item $\brs{\phi} \leq A$.
\end{enumerate}
Then there exists $C = C(h, A)$ so that, setting $dm = e^{- \phi} dV_g$, every smooth real function $u$ satisfies
\begin{align} \label{f:weightedSobolev}
\left( \int_M \brs{u}^{4n/(2n-1)} dm \right)^{(2n-1)/(2n)} \leq C \int_M \left( \brs{\del v}^2_g + v^2 \right) dm.
\end{align}
\begin{proof} We first observe using the bound on $\phi$ and elementary linear algebra that
\begin{align} \label{l:sob10}
\brs{\del u}_h^2 dV_h \leq C (\tr_g h)^{n-1} \brs{\del u}_g^2 dm, \qquad dV_h \leq C (\tr_g h)^n dm.
\end{align}
Now let $s = \frac{8n}{4n+1} < 2$, and consider H\"older exponents $p = 2/s, q = (2/(2-s))$.  Using (\ref{l:sob10}) and H\"older's inequality we obtain
\begin{align*}
    \int_M \brs{\del u}^{s}_h dV_h \leq&\ C \int_M \left( \brs{\del u}^2_g e^{-\phi} \frac{\gw^n}{\gw_h^n} \right)^{s/2} (\tr_g h)^{(n-1) s/2} dV_h\\
    \leq&\ C \left( \int_M \brs{\del u}^2_g dm \right)^{s/2} \left( \int_M (\tr_g h)^{(n-1) s/(2-s)} dV_h \right)^{(2-s)/2}.
\end{align*}
Taking the $2/s$ power of this and applying hypothesis (1) yields
\begin{align} \label{l:sob20}
\brs{\brs{\del u}}_{L^s(h)}^2 \leq C \left( \int_M (\tr_g h)^{4n(n-1)} dV_h \right)^{1/4n} \int_M \brs{\del u}^2 dm \leq C \int_M \brs{\del u}^2 dm.
\end{align}
A similar argument yields
\begin{align} \label{l:sob30}
\brs{\brs{u}}_{L^s(h)}^2 \leq C \left( \int_M (\tr_g h)^{4 n^2} dV_h \right)^{1/(4n)} \int_M u^2 dm.
\end{align}
Applying the Sobolev inequality $W^{1,s} \to L^{8n/(4n-3)}$ for $h$ and using (\ref{l:sob20}) and (\ref{l:sob30}) thus yields
\begin{align} \label{l:sob40}
    \brs{\brs{u}}_{L^{8n/(4n-3)}(h)}^2 \leq C \int_M \left( \brs{\del u}^2_g + u^2 \right) dm.
\end{align}
Finally we estimate using hypotheses (1) and (2) and H\"older's inequality,
\begin{align*}
    \int_M \brs{u}^{4n/(2n-1)}dm =&\ \int_M \brs{u}^{4n/(2n-1)} e^{-\phi} \frac{\gw^n}{\gw_h^n} dV_h\\
    \leq&\ \left( \int_M \left( e^{-\phi} \frac{\gw^n}{\gw_h^n} \right)^{4n-2} dV_h \right)^{1/(4n-2)} \left( \int_M \brs{u}^{8n/(4n-3)} dV_h \right)^{(4n-3)/(4n-2)}\\
    \leq&\ C \brs{\brs{u}}_{L^{8n/(4n-3)}(h)}^{4n/(2n-1)}.
\end{align*}
Taking the $(2n-1)/2n$-power of this inequality and returning to (\ref{l:sob40}) finishes the proof.
\end{proof}
\end{lemma}

Next we adapt the classic parabolic Moser estimate \cite{moser1964harnack} to our setting, using a certain weighted measure adapted to the scalar curvature lower bound, and moreover using the weighted Sobolev inequality of Lemma \ref{l:Sobolev}.

\begin{lemma} \label{l:moser} Fix $(M, h,J)$ a compact Hermitian manifold, and suppose $(g_t, \phi_t)$ is a one-parameter family of Hermitian metrics and smooth functions on $[0, \tau], \tau \leq 1$, which satisfy the hypotheses of Lemma \ref{l:Sobolev} uniformly for all $t$.  Assume that for every nonnegative smooth function $u$ one has
\begin{align} \label{f:moser10}
    \frac{d}{dt} \int_M u dm \leq \int_M \left( \square u + K u \right) dm, \qquad K \geq 0.
\end{align}
If $v > 0$ satisfies $\square v \leq A (1 + \tr_g h) v$, then there exists $C = C(h, A)$ such that
\begin{align} \label{f:moser20}
    \sup_{M \times [\tau/2, \tau]} v \leq [C (p^2 + K + \tau^{-1})]^{(2n+1)/p} \left( \int_0^\tau \int_M v^p dm dt \right)^{1/p}.
\end{align}
\begin{proof} We first observe the pointwise identity
\begin{align*}
    \square v^p = p v^{p-1} \square v - p(p-1) v^{p-2} \brs{\del v}^2_g = p v^{p-1} \square v - \frac{4(p-1)}{p} \brs{\del v^{p/2}}^2_g.
\end{align*}
Using this, (\ref{f:moser10}), and our assumption on $v$,
\begin{align} \label{f:moser30}
\frac{d}{dt} \int_M v^p dm + \frac{4(p-1)}{p} \int_M \brs{\del v^{p/2}}^2_g dm \leq&\ \int_M \left[ Cp(1 + \tr_g h) + K \right] v^p dm.
\end{align}
By Lemma \ref{l:Sobolev}, H\"older's inequality and interpolation, we observe that every smooth function $w$ satisfies
\begin{align*}
\int_M (\tr_g h) w^2 dm \leq&\ C \brs{\brs{w}}_{L^{8n/(4n-1)}(dm)}^2\\
\leq&\ C \brs{\brs{w}}_{L^{4n/(2n-1)}(dm)} \brs{\brs{w}}_{L^2(dm)}\\
\leq&\ \ge \int_M \brs{\del w}^2_g dm + C \ge^{-1} \int_M w^2 dm.
\end{align*}
Choosing $w = v^{p/2}$ and choosing $\ge = \delta p^{-1}$ for small $\delta$ then yields
\begin{align} \label{f:moser40}
    \frac{d}{dt} \int_M v^p dm + \int_M \brs{\del v^{p/2}}^2_g dm \leq C \left(p^2 + K \right) \int_M v^p dm.
\end{align}
We also observe via interpolation and Lemma \ref{l:Sobolev} that for any smooth positive function $w$ we have for any time interval $I$,
\begin{align} \label{f:moser50}
    \int_I \int_M w^{2 + 1/n} dm dt \leq&\ C \left( \sup_I \int_M w^2 dm \right)^{1/(2n)} \int_I \int_M \left( \brs{\del w}^2_g + w^2 \right) dm dt.
\end{align}
Combining (\ref{f:moser40}) and (\ref{f:moser50}) allows for the standard iteration with exponents $p_j = p(1 + 1/(2n))^j$ and time intervals $[(\tau/2)(1 - 2^{-j}), \tau]$, yielding the claimed estimate.
\end{proof}
\end{lemma}

\begin{proof}[Proof of Theorem \ref{t:mainthm1}] We fix any time $\tau < \tau^*$ and construct a solution $(\beta, f)$ against a uniformly controlled moving background on $[0, \tau]$, and establish uniform estimates for $g$ on this time interval.  By Corollary \ref{c:Lptrace} we obtain a uniform $L^p$ bound on $\tr_g h$ for any $p$.  Moreover, by Proposition \ref{p:dilatonscalarbound} we have uniform upper and lower bounds on the dilaton, so we have justified the hypotheses of Lemma \ref{l:Sobolev} uniformly for all $t$.  Furthermore, using the a priori lower bound on the weighted scalar curvature $R^{H,\phi}$ from Proposition \ref{p:dilatonscalarbound} and equation (\ref{f:phimeasureevolution}), we have justified the hypothesis (\ref{f:moser10}) of Lemma \ref{l:moser}.  Since $\tr_g h$ is a valid test function for Lemma \ref{l:moser} by the Schwarz Lemma \cite[Lemma 6.2]{StreetsPCFBI}, we obtain a uniform upper bound for $\tr_g h$.  As explained in the proof of Proposition \ref{p:fmean} we also have an a priori upper bound for the volume form, and hence uniform parabolicity.  Higher order estimates for $g$ follow from \cite[Theorem 1.7]{StreetsPCFBI}.
\end{proof}

\section{Proof of Theorem \ref{t:AubinYau}}
\subsection{Normalized reduction and maximum principles}

Our first goal is to obtain a refined $(\beta, f)$-reduction for the normalized flow.  Fix $(M^{2n}, J)$ a compact K\"ahler manifold with $c_1(M) < 0$, and fix a pluriclosed metric $\gw_0$.  By Theorem \ref{t:mainthm1}, since $[\omega_0] + t[\gw_h]$ is a positive class for all $t \geq 0$, we know that the solution $\til{\gw}_t$ to pluriclosed flow with initial condition $\gw_0$ exists smoothly for $t \geq 0$.  Now consider $s = \log (1 + t)$, and let $\gw(s) = e^{-s} \til{\gw}(e^s - 1)$.  It follows that
\begin{align*}
    \frac{\del}{\del s} \gw = - \rho_B^{1,1}(\omega) - \gw,
\end{align*}
that is, $\gw(s)$ is a solution of the normalized pluriclosed flow.

We now construct a $(\beta, f)$ reduction for the normalized flow following \S \ref{ss:reduction}.  Choose a background K\"ahler metric $h$ with $\gw_h \in - c_1(M)$.  By the $\del\delb$-lemma there exists a smooth function $\chi$ such that $d d^c \chi = \rho(h) + \omega_h$.  Let $\nu = e^{\chi} \gw_h^n$, with $\chi$ normalized so that $\nu$ has unit volume, and then observe $\rho(\nu) = - \gw_h$.  Furthermore, by Lemma \ref{l:cones} we fix a K\"ahler metric $\gw_K \in [\gw_0]$, and define the moving K\"ahler background $\til{\gw}_t := \gw_K + t \gw_h$.  We further define $\chi(t) = t \chi$, and obtain $(\til{\beta}, \til{f})$ a solution to (\ref{f:reducedflow}).  Define
\begin{align*}
    \hat{\gw}(s) = e^{-s} \til{\gw}(e^s - 1), \qquad \gb(s) = e^{-s} \til{\gb}(e^s - 1), \qquad f(s) = e^{-s} \til{f}(e^s - 1) - n(s - 1 + e^{-s}).
\end{align*}
A straightforward computation shows that
\begin{gather} \label{f:normalizedreduction}
    \begin{split}
        \hat{\gw}_s =&\ e^{-s} \gw_K + (1 - e^{-s}) \gw_h,\\
        \square \gb =&\ - T_{g_s} \circ \delb \beta - \gb,\\
        \frac{\del}{\del s} f =&\ \log \frac{\gw^n}{\nu} - \tr_{\gw_s} (\del \bgb + \delb \gb) - f,
    \end{split}
\end{gather}
where here and throughout this section
\begin{align*}
    \square = \frac{\del}{\del s} - \Delta^C_{g_s}
\end{align*}
is the heat operator associated to the time dependent Chern Laplacian associated to $g_s$.

\begin{prop} \label{p:AYflower} Given the setup above, one has.
\begin{align*} 
    \square f =&\ \tr_{g_s} \hat{g} + \log \frac{\gw_s^n}{\nu} - n - f.
\end{align*}
Moreover, for all $s \geq 0$ we have
\begin{align*}
    \inf_{M \times \{s\}} f \geq - C.
\end{align*}
\begin{proof} The claimed heat equation follows from (\ref{f:normalizedreduction}).  Arguing as in Proposition \ref{p:flower}, if we let $\gl_i$ denote the eigenvalues of $g^{-1} \hat{g}$ we then obtain
\begin{align*}
    \square f = \sum_i \left(\gl_i - 1 - \log \gl_i \right) + \log \frac{\hat{\gw}^n}{\gw_h^n} - f \geq - C - f,
\end{align*} 
and the estimate follows by the maximum principle.
\end{proof}
\end{prop}

\begin{prop} \label{p:AYbetabounds} Given the setup above, one has
\begin{align*}
    \square \brs{\beta}^2_g =&\ - \brs{\N \beta}^2_g - \brs{\bar{\N} \beta + \beta \cdot \bar{T}}^2_g - \brs{\gb}^2_g\\
    \square \brs{\del \beta}^2_g =&\ - \brs{\N \del \beta}^2_g - \tfrac{1}{2} \brs{T}^2_g - \IP{Q, \tr_g \del \beta \otimes \delb \bar{\beta}}_g.
\end{align*}
Moreover, for all $s \geq 0$ we have
\begin{align*}
    \sup_{M \times \{s\}} \brs{\beta}^2_g + \brs{\del \beta}^2_g \leq&\ C e^{-s}.
\end{align*}
\begin{proof}
The evolution equations follow from those in Proposition \ref{p:betabounds} by rescaling.  Also exponential decay of $\brs{\beta}^2_g$ follows by Proposition \ref{p:betabounds} and the scaling law.  To establish exponential decay of $\brs{\del \beta}^2_g$ we exploit the $\brs{\N \beta}^2_g$ term in the evolution of $\brs{\gb}^2$.  By Cauchy-Schwarz and interpolation it follows that
\begin{align*}
    \brs{\del \beta}^2_g \leq 4 \brs{\N \gb}^2_g + \brs{\gb}^2_g \brs{T}^2_g.
\end{align*}
Using exponential decay of $\brs{\beta}^2_g$ it follows that for sufficiently large times $s$ we have
\begin{align*}
    \square \left( \brs{\del \beta}^2 + 4 \brs{\gb}^2 \right) \leq&\ - \tfrac{1}{2} \brs{T}^2_g - 4 \brs{\N \beta}^2_g - 4 \brs{\gb}^2_g \leq \left( - \brs{\del \beta}^2_g + 4 \brs{\gb}^2_g \right).
\end{align*}
The result follows from the maximum principle.
\end{proof}
\end{prop}

\begin{prop} \label{p:AYvolumelower} Given the setup above,
\begin{align*}
    \square \log \frac{\gw^n}{\nu}  = \brs{T}^2_g + \tr_g h - n.
\end{align*}
Moreover, for all $s \geq 0$ we have
\begin{align*}
    \inf_{M \times \{s\}} \log \frac{\gw^n}{\nu} \geq - C.
\end{align*}
\begin{proof} The claimed evolution equation is standard computation (cf. \cite[Lemma 6.1]{StreetsPCFBI}).  Given this and the inequality
\begin{align*}
    \tr_g h \geq n \left( \frac{\gw_h^n}{\nu} \right)^{1/n} e^{- \frac{1}{n} \log (\gw^n/\nu)},
\end{align*}
the claimed lower bound follows from the maximum principle.
\end{proof}
\end{prop}

\begin{prop} \label{p:AYSchwarz} Given the setup above, we have
\begin{align} \label{f:AYSchwarz10}
\tr_g h \leq C e^{Cf}, \qquad e^{2 \brs{\del \beta}^2_g} \frac{\gw^n}{\nu} \leq C e^{C f}.
\end{align}
\begin{proof} Taking the normalization into account in the estimate of Corollary \ref{c:Lptrace} we obtain
\begin{align*}
\square \tr_g h \leq C (\tr_g h)^2 + \tr_g h, \qquad \square \log \tr_g h \leq C \tr_g h + 1.
\end{align*}
Hence for large $A$ we estimate using Proposition \ref{p:AYflower} and the lower bound of Proposition \ref{p:AYvolumelower},
\begin{align*}
    \square \left( \log \tr_g h - A f \right) \leq&\ - \left(A c - C \right) \tr_g h + 1 + An - A \log \frac{\gw^n}{\nu} + A f\\
    \leq&\ C - \left( \log \tr_g h - A f \right)
\end{align*}
The first claimed inequality follows by the maximum principle.
Similarly, using Propositions \ref{p:AYflower}, \ref{p:AYbetabounds}, and \ref{p:AYvolumelower} we estimate for large $A$,
\begin{align*}
    \square \left( \log \frac{\gw^n}{\nu} + 2 \brs{\del \beta}^2_g - A f \right) \leq&\ - (A c - 1) \tr_g h + n(A - 1) - (A - 1) \log \frac{\gw^n}{\nu} + 2 \brs{\del \beta}^2_g\\
    &\ - \left( \log \frac{\gw^n}{\nu} + 2 \brs{\del \beta}^2 - A f \right)\\
    \leq&\ C - \left( \log \frac{\gw^n}{\nu} + 2 \brs{\del \beta}^2 - A f \right).
\end{align*}
The second claimed inequality follows from the maximum principle.  
\end{proof}
\end{prop}

\subsection{Potential mean estimate}

Next we obtain an estimate for the average value of $f$.  The proof exploits a monotonicity for a relative entropy functional in a similar spirit to \cite{berman2013relative}.

\begin{prop} \label{p:AYmeanest} Given the setup above,
\begin{align*}
    \sup_{s} \brs{ \bar{f}_h} \leq C. 
\end{align*}
\begin{proof} Let
\begin{align*}
    J(s) := \int_M \log \frac{\gw^n}{\nu} d \nu.
\end{align*}
We observe using Jensen's inequality and Lemma \ref{l:volumebound} that
\begin{align} \label{f:AYmean5}
    J \leq \log \int_M \frac{\gw^n}{\nu} d \nu = \log \int_M \gw^n = C.
\end{align}
Using Proposition \ref{p:AYvolumelower} and an integration by parts argument as in Lemma \ref{l:weightbounds}, we obtain
\begin{align} \label{f:AYmean10}
    \frac{d}{ds} J = \int_M \left( \brs{\theta - \del \log \frac{\gw^n}{\nu}}^2_g + \tfrac{1}{2} \brs{T}^2 + \tr_g h \right) d \nu - n.
\end{align}
Using the evolution equation (\ref{f:normalizedreduction}) for $f$ and integrating by parts gives
\begin{gather} \label{f:AYmean20}
    \begin{split}
        \frac{d}{ds} \int_M f d \nu =&\ - \int_M f d \nu + J - \int_M \tr_{\gw} \left( \del \bgb + \delb \gb \right) d \nu\\
        =&\ - \int_M f d \nu + J + 2 \int_M \mathrm{Im} \IP{ \gb, \theta - \del \log \frac{\gw^n}{\nu}}_g d \nu\\
        \leq&\ - \int_M f d \nu + J + \int_M \left( \brs{\gb}^2_g + \brs{\theta - \del \log \frac{\gw^n}{\nu}}^2_g \right) d \nu.
    \end{split}
\end{gather}
Combining (\ref{f:AYmean10}) and (\ref{f:AYmean20}), and using Proposition \ref{p:AYbetabounds} we obtain
\begin{align*}
    \frac{d}{ds} \left( \int_M f d \nu - J \right) \leq&\ n + \int_M \brs{\gb}^2_g d \nu - \left(\int_M f d \nu - J \right) \leq n + C e^{-s} - \left(\int_M f d \nu - J \right).
\end{align*}
Since $J$ is bounded above by (\ref{f:AYmean5}), it follows that $\int_M f d \nu$ is uniformly bounded above.  Since $f$ is bounded below by Proposition \ref{p:AYflower} the proposition follows easily.
\end{proof}
\end{prop}

\subsection{Main Proof}

\begin{proof}[Proof of Theorem \ref{t:AubinYau}] The main step is to establish uniform parabolicity of the normalized flow via Moser iteration as in the proof of Theorem \ref{t:mainthm1}.  A further complication here is that the scalar curvature lower bound is a priori  degenerating due to the scaling.  We choose $\til{\phi}_t$ a solution to (\ref{f:dilatonflow}) and then let $\phi(s) = \til{\phi}(e^s - 1)$, and set $dm = e^{-\phi} dV_g$.  By Proposition \ref{p:dilatonscalarbound} and rescaling we get
\begin{align} \label{f:AYpf10}
\square \phi = \tfrac{1}{2} \brs{T}^2_g, \qquad 0 \leq \phi \leq C, \qquad R^{H,\phi} \geq - C e^s.
\end{align}
By incorporating the scaling we see that for an arbitrary smooth function $u$ we have
\begin{align} \label{f:AYpf20}
    \frac{d}{ds} \int_M u dm =&\ \int_M \left( \square u - \left( \tfrac{1}{2} R^{H,\phi} + n\right) u \right) dm.
\end{align}
In particular, due to (\ref{f:AYpf10}) we cannot verify the bound of Lemma \ref{l:moser} in a uniform way for all times.  We rectify this by starting the iteration at values of $p$ which grow with time.

We first prove a uniform $C^0$ bound for $\tr_g h$.  We first observe using Proposition \ref{p:AYSchwarz}, and applying the Gaussian bound of Proposition \ref{p:Gaussian}, and the estimates of Propositions \ref{p:AYbetabounds} and \ref{p:AYmeanest},
\begin{align*}
    \int_M (\tr_g h)^p \leq C^{p+1} \int_M e^{C(p+1) f} dV_h \leq&\ \exp \left( C p( 1 + \bar{f}_h + p \sup_M \brs{\gb}^2_g) \right) \leq \exp( Cp + C p^2 e^{-s}).
\end{align*}
Notice that this bound, and a similar bound for $\log \frac{\gw^n}{\nu}$, implies that the Sobolev inequality of Lemma \ref{l:Sobolev} can be applied uniformly at all times.  Now consider the time interval $[\sigma, \sigma + 1]$ for large $\sigma$.  Using the evolution equation for $\tr_g h$ from Proposition \ref{p:AYSchwarz}, we apply Lemma \ref{l:moser} with $p = \sigma + 2$ to yield
\begin{align*}
    \sup_{M \times [\sigma + \tfrac{1}{2}, \sigma + 1]} \tr_g h \leq&\ \left[ C (p^2 + e^{\sigma + 1} + 1 \right)]^{(2n+1)/p} \exp \left(C + C p e^{-\sigma} \right) \leq C.
\end{align*}
An identical argument applies to yield a uniform upper bound for $e^{2 \brs{\del \beta}^2_g} \frac{\gw^n}{\nu}$, which implies uniform parabolicity of the flow.  Given this the higher regularity theory of \cite{StreetsPCFBI} applies to yield uniform a priori $C^k$ estimates for $g$.  By a standard argument using the expander entropy functional for generalized Ricci flow \cite{Streetsexpent}, we obtain $s_i \to \infty$ such that $\{g_{s_i}\}$ converges to a K\"ahler-Einstein metric $g_{\infty}$, which is the unique such in $- c_1(M)$ by the maximum principle \cite{calabi1957kahler}.  By \cite[Theorem 1.2]{HCF} the flow converges exponentially as claimed.
\end{proof}

\section{Proof of Theorem \ref{t:CalabiYau}}
\subsection{Refined reduction and maximum principles}

Fix $(M^{2n}, J)$ a compact K\"ahler manifold with $c_1(M) = 0$, and fix a pluriclosed metric $\gw_0$.  We first refine the $(\beta, f)$ reduction in this setting.  Using Lemma \ref{l:cones} choose a background K\"ahler metric $h$ with $\gw_h \in [\gw_0]$.  Since $c_1(M) = 0$, by the $\del\delb$-lemma there exists a smooth function $\chi$ such that $d d^c \chi = \rho(h)$.  Without loss we can assume $\sup \chi = 0$, and define $\nu = e^{\chi} \gw_h^n$, hence $\rho(\nu) = 0$.  Following the notation of \S \ref{s:Kahlercone}, we define $\hat{\gw}(t) = \gw_h$, $\chi(t) = t \chi$, resulting in the system
\begin{gather} \label{f:CYsystem}
    \begin{split}
        \square \gb =&\ - T_{g_t} \circ \delb \gb,\\
        \dt f =&\ \log \frac{\gw_t^n}{\nu} - \tr_{\gw_t}(\del \bgb + \delb \gb).
    \end{split}
\end{gather}

\begin{prop} \label{p:CYMP} Given the setup above,
\begin{align*}
    \sup_{M \times \{t\}} \brs{\log \frac{\gw^n}{\gw_h^n}} \leq C.
\end{align*}
\begin{proof} Using that $\rho(\nu) = 0$, we have
\begin{align*}
    \square \log \frac{\gw^n}{\nu} = \brs{T}^2_g \geq 0.
\end{align*}
Hence by the maximum principle there is an a priori lower bound for $\log \frac{\gw^n}{\nu}$, hence also for $\log \frac{\gw^n}{\gw_h^n}$.  Similarly using Proposition \ref{p:betabounds} we have
\begin{align*}
    \square \left( \log \frac{\gw^n}{\nu} + 2 \brs{\del \beta}^2_g \right) \leq 0,
\end{align*}
which yields the claimed upper bound by the maximum principle.
\end{proof}
\end{prop}

\subsection{Extension of Yau \texorpdfstring{$C^0$}{C0} estimate}

Here we extend Yau's $C^0$ estimate for the complex Monge-Amp\`ere equation \cite{YauCC} to our setting coupled to $\beta$.  In particular we apply a Moser iteration argument taking care of the extra $\beta$ terms, to yield a lower bound on the normalized potential $f - \bar{f}_h$.  We begin with a linear algebra lemma.

\begin{lemma} \label{l:YauC0LA} Fix $F$ a real two-form such that $F^{1,1} = \gw > 0$.  Given $\gw_h > 0$ another K\"ahler form and $\ga$ a real $1$-form define
\begin{align*}
    W_k = F^k \wedge \gw_h^{n-1-k}, \qquad V_k = F^k \wedge \gw_h^{n-k}, \qquad E_k(\ga) = \tfrac{1}{2} \ga \wedge J \ga \wedge W_k.
\end{align*}
Given $\gb$ a $(1,0)$-form such that $\brs{\gb}^2_g \leq B$, one has
\begin{align*}
    \brs{\ga \wedge (\beta + \bar{\gb}) \wedge W_k}^2 \leq \frac{4 B}{k+1} E_k(\ga) V_{k+1},
\end{align*}
interpreted via densities relative to any fixed volume form.
\begin{proof} We first observe by type considerations that $V_k > 0$ and $E_k(\ga) \geq 0$.  Now define
\begin{align*}
    F_s = F + s \ga \wedge (\gb + \bgb) + \tfrac{1}{2} B s^2 \ga \wedge J \ga.
\end{align*}
It follows easily from $\brs{\gb}^2_g \leq B$ that $F^{1,1}_s \geq 0$.  Using that $(F_s - F)^2 = 0$ we compute
\begin{align*}
    0 \leq F_s^{k+1} \wedge \gw_h^{n-1-k} = V_{k+1} + (k+1) s \ga \wedge (\beta + \bgb) \wedge W_k + (k + 1) B s^2 E_k(\ga).
\end{align*}
Interpreting the above inequality as a quadratic polynomial in $s$, it follows that the discriminant is nonpositive, giving the claim.
\end{proof}
\end{lemma}

\begin{prop} \label{p:CYpotentialLB} Fix $(M^{2n}, J, h)$ a compact K\"ahler manifold and let
\begin{align*}
F = \gw_h + d (\beta + \bgb + d^c f),
\end{align*}
where $\gb \in \Lambda^{1,0}$ and $F^{1,1} = \gw > 0$.  Assume
\begin{enumerate}
    \item $F^n \leq A \gw_h^n$,
    \item $\brs{\gb}^2_g \leq B$.
\end{enumerate}
Then there exists $C = C(n, h, A, B)$ such that $f - \bar{f}_h \geq - C$.
\begin{proof} We set $v = 1 + (\bar{f}_h - f)_+ \geq 1$.  Fix $p \geq 1$ and let
\begin{align*}
    \mathcal I_k(p) = \int_M v^p F^k \wedge \gw_h^{n-k}, \qquad \mathcal E_k(p) = \int_{f < \bar{f}_h} v^{p-1} df \wedge d^c f \wedge F^k \wedge \gw_h^{n-1-k}.
\end{align*}
First observe by Stokes Theorem and Lemma \ref{l:YauC0LA} (using hypothesis (2)) that
\begin{align*}
    \mathcal I_{k+1}(p) - \mathcal I_k(p) =&\ \int_M v^p d \left( \gb + \bgb + d^c f \right) \wedge F^k \wedge \gw_h^{n-1-k}\\
    =&\ p \int_{\{f < \bar{f}_h\}}  v^{p-1} df \wedge (\gb + \bgb + d^c f) \wedge F^k \wedge \gw_h^{n-1-k}\\
    \geq&\ \frac{p}{2} \mathcal E_k(p) - \frac{2 B p}{k+1} \mathcal I_{k+1}.
\end{align*}
This yields in particular an upper bound for $\mathcal E_0(p)$, and then an induction argument together with hypothesis (1) yields
\begin{align*}
    \mathcal E_0(p) \leq \frac{c_n (1 + 2 B p)^n}{p} \mathcal I_n(p) \leq \frac{c_n A (1 + 2 B p)^n}{p} \mathcal I_0(p).
\end{align*}
Expressing this in terms of $v$ and using $dv = -df$ on $\{f < \bar{f}_h\}$ and $v \geq 1$ we obtain
\begin{align*}
    \int_M \brs{d v^{(p+1)/2}}^2_{h} dV_h \leq C (p + 1)^{n+1} \int_M v^{p+1} dV_h.
\end{align*}
A standard Moser iteration argument and Proposition \ref{p:Gaussian} then yields $\brs{\brs{v}}_{\infty} \leq C \brs{\brs{v}}_{L^1(h)} \leq C$.
\end{proof}
\end{prop}

\subsection{Main Proof} 

A final lemma necessary for the main proof is a priori control over the variation of $\bar{f}_h$ over bounded intervals.  The argument is similar in spirit to Proposition \ref{p:fmean}.

\begin{lemma} \label{l:CYaveragebnds} Given the setup above, for $\brs{t_1 - t_2} \leq 1$ we have
\begin{align*}
    \brs{\bar{f}_h(t_1) - \bar{f}_h(t_2)} \leq C.
\end{align*}
\begin{proof} We first obtain a lower bound for $f$.  Arguing as in Proposition \ref{p:flower}, we let $\gl_i$ denote the eigenvalues of $g h^{-1}$ and use $\nu \leq \gw_h^n$ to obtain
\begin{align} \label{f:CYfheat}
    \square f =&\ \tr_g h + \log \frac{\gw^n}{\nu} - n =  \sum_i \left( \gl_i^{-1} + \log \gl_i - 1 \right) + \log \frac{\gw_h^n}{\nu} \geq 0.
\end{align}
Hence the minimum of $f$ is nondecreasing.  Combining this with Proposition \ref{p:CYpotentialLB} yields an a priori lower bound on $\bar{f}_h(t_2) - \bar{f}_h(t_1)$.

To obtain the upper bound we first exploit the lower bound on $f$ just shown together with Proposition \ref{p:CYpotentialLB} and define
\begin{align*}
    \til{f} = f - \bar{f}_h(t_1) + C + 1 \geq 1, \qquad A(t) = \int_M \til{f} d \mu,
\end{align*}
where $d \mu = e^{-2 \phi} dV_g$ (cf. Lemma \ref{l:weightbounds}).  We first observe by Propositions \ref{p:dilatonscalarbound} and \ref{p:Gaussian} that $A(t_1) \leq C$.  Using (\ref{f:CYsystem}) and Lemma \ref{l:weightbounds}, and estimating as in Proposition \ref{p:fmean} we obtain
\begin{align*}
    \frac{d}{dt} A \leq K A + C + C \int_M \left( \brs{T}^2_g + \brs{\del \phi}^2_g \right) d \mu.
\end{align*}
Since the mass of $\mu$ is uniformly controlled above and below the argument of Lemma \ref{l:weightbounds} yields an a priori bound on the time integral of the final term above over an interval of unit length.  This yields an a priori upper bound for $A$, and finally exploiting the $L^1$ bound on $f - \bar{f}_h$ from Proposition \ref{p:Gaussian} gives the upper bound for $\bar{f}_h(t_2) - \bar{f}_h(t_1)$.
\end{proof}
\end{lemma}

\begin{proof} [Proof of Theorem \ref{t:CalabiYau}] The main goal is to obtain a uniform $L^p$ estimate for $\tr_g h$, after which the proof follows the Moser iteration argument of Theorem \ref{t:mainthm1}.  The key new point is to apply the Schwarz Lemma with the potential normalized by a time-average.  Fix a large constant $A$, $t \geq 1$, and let
\begin{align*}
    \Psi = \log \tr_g h - A \left( f - \int_{t-1}^t \bar{f}_h(s) ds \right).
\end{align*}
Using the Schwarz Lemma as in Corollary \ref{c:Lptrace}, equation (\ref{f:CYfheat}), and applying Proposition \ref{p:CYMP} and Lemma \ref{l:CYaveragebnds} we obtain
\begin{align*}
    \square \Psi \leq&\ \left(C - A \right) \tr_g h + A \left(n - \log \frac{\gw^n}{\nu} \right) + A \left( \bar{f}_h(t) - \bar{f}_h(t-1) \right) \leq - \tr_g h + C.
\end{align*}
Using the maximum principle and Proposition \ref{p:CYpotentialLB} we obtain $\Psi \leq C$ uniformly for all times $t$.  Using Lemma \ref{l:CYaveragebnds} we may replace the spacetime average term of $f$ by $\bar{f}_h(t)$ to obtain the uniform a priori estimate
\begin{align*}
    \tr_g h \leq C e^{A(f - \bar{f}_h)}.
\end{align*}
By Proposition \ref{p:Gaussian} we hence obtain a uniform $L^p$ estimate for $\tr_g h$ for all times.

We have thus established the required estimates to repeat the Moser iteration argument in the proof of Theorem \ref{t:mainthm1}, yielding a uniform $C^0$ estimate for $\tr_g h$, which then by Proposition \ref{p:CYMP} implies uniform upper and lower bounds for $g$ at all times.  By the higher regularity theory of \cite{StreetsPCFBI} we obtain uniform estimates of all $C^k$ norms of $g$.  By a standard argument using the $\gl$-functional monotonicity for generalized Ricci flow (\cite{OSW}, \cite[Chapter 6]{GRFbook}) we obtain $t_i \to \infty$ such that $\{g(t_i)\}$ converges to $g_{\infty}$ a pluriclosed Ricci soliton.  Observing that the closed Bismut torsion of a pluriclosed metric is exact on a K\"ahler background by the $\del\delb$-lemma, it follows from \cite[Proposition 2.1]{streets2026topology} that the torsion in fact vanishes, and $g_{\infty}$ is K\"ahler, Ricci-flat.  It is the unique such metric in the K\"ahler class by \cite{Calabiaffine}.  With this subsequential convergence established exponential convergence of the flow follows from \cite[Theorem 1.2]{HCF}
\end{proof}

\subsection{Contractibility of generalized K\"ahler classes}

\begin{proof}[Proof of Corollary \ref{c:GKcor}] We refer to \cite{AFSU} for relevant background on generalized K\"ahler structures and their relationship to Calabi-Yau manifolds.  Fix $(M^{2n}, J)$ a compact K\"ahler Calabi-Yau manifold, $\sigma$ a holomorphic Poisson tensor, and admissible spinor classes $\ga, \gb$ such that $\mathcal {GK}^{J, \sigma}_{\ga,\gb}$ is nonempty.  Fix $(g_0, I_0, J) \in \mathcal {GK}^{J, \sigma}_{\ga,\gb}$, and let $(g_t, I_t, J)$ denote the unique solution to generalized K\"ahler-Ricci flow \cite{GKRF} in the $J$-fixed gauge.  Since $\gw_J$ is in particular a solution to pluriclosed flow, by Theorem \ref{t:CalabiYau} $(g_t, I_t)$ exists for all time and $\gw_J$ converges to a unique K\"ahler Calabi-Yau metric, and furthermore we can gauge-fix the flow of the $b$-field so that the limit $b_{\infty} = 0$.  This limiting metric, and the entire generalized K\"ahler structure, is canonically determined by the spinor classes $\ga, \gb$ by \cite[Proposition 2.14]{AFSU}, finishing the proof.
\end{proof}


\end{document}